\documentclass[12pt,draft]{amsart}
\usepackage{amsmath,amsthm,latexsym,amscd,amsbsy,amssymb,fleqn,leqno}
\chardef\bslash=`\\ 

\makeatletter
\def\verbatim{\interlinepenalty\@M \@verbatim
  \leftskip\@totalleftmargin\advance\leftskip2pc
  \frenchspacing\@vobeyspaces \@xverbatim}
\makeatother
\makeatletter
  \def\dgt@k{\dg@DX=-3 \dg@DY=2 \dg@SIZE=3} 
\makeatother

\makeatletter
  \def\dgt@kk{\dg@DX=3 \dg@DY=-1 \dg@SIZE=3}%
\makeatother

\theoremstyle{plain}
\newtheorem{thm}{Theorem}[section]

\newtheorem{lem}[thm]{Lemma}

\theoremstyle{definition}

\numberwithin{equation}{section}

\newcounter{rmnum}

\def\symbolnote#1#2{\let\thefootn=\thefootnote%
\renewcommand{\thefootnote}{\fnsymbol{footnote}}%
\footnotemark[#1]%
\footnotetext[#1]{#2}%
\let\thefootnote=\thefootn
}

\newfont{\bbb}{msbm10 scaled \magstep1}
\newfont{\bbc}{msbm8 scaled \magstep0}

\newcommand{\N}{\mbox{\bbb N}}

\newcommand{\uin}{\mbox{\bbb I}}
\newcommand{\e}{\mbox{\rm e-dim}}

\begin{document}


\title{Extensional dimension and completion of maps}

\author{H. Murat Tuncali}
\address{Department of Computer Science and Mathematics,
Nipissing University,
100 College Drive, P.O. Box 5002, North Bay, ON, P1B 8L7, Canada}
\email{muratt@nipissingu.ca}
\thanks{The authors were partially supported by their NSERC grants.}

\author{E. D. Tymchatyn}
\address{Department of Mathematics and Statistics, University of Saskatchewan,
McLean Hall, 106 Wiggins Road, Saskatoon, SK, S7N 5E6, Canada}
\email{tymchat@math.usask.ca}

\author{Vesko Valov}
\address{Department of Computer Science and Mathematics, Nipissing University,
100 College Drive, P.O. Box 5002, North Bay, ON, P1B 8L7, Canada}
\email{veskov@nipissingu.ca}

 

\begin{abstract}
We prove the following completion theorem for closed maps between metrizable spaces: 
Let $f\colon X\to Y$ be a closed surjection between metrizable spaces with $\e f\leq K$, $\e X\leq L_X$ and $\e Y\leq L_Y$  for some countable $CW$-complexes $K$, $L_X$ and $L_Y$. Then there exist completions $\widetilde{X}$ and $\widetilde{Y}$ of $X$ and $Y$, respectively, and a closed surjection $\widetilde{f}\colon\widetilde{X}\to\widetilde{Y}$ extending $f$ such that $\e\widetilde{f}\leq K$, $\e\widetilde{X}\leq L_X$ and $\e\widetilde{Y}\leq L_Y$. We   
also establish a parametric version of a result of Katetov characterizing the covering dimension of metrizable spaces in terms of uniformly 0-dimensional maps into finite-dimensional cubes.  
\end{abstract}

\maketitle

\markboth{H. M.~Tuncali, E. Tymchatyn and V.~Valov}{Completion of maps}


\section{Introduction}

Katetov \cite{mk} and Morita \cite{km} proved that every finite-dimensional metrizable space has a metrizable completion of the same dimension. A completion theorem for extensional dimension with respect to countable $CW$-complexes was established by Olszewski \cite{wo} in the class of separable metrizable spaces and recently by Levin \cite{ml1} in the class of all metrizable spaces.    

Concerning completions of maps with the same dimension, Keesling \cite{jk} proved that if $f\colon X\to Y$ is a closed surjective map between metrizable finite-dimensional spaces, then there are completions $\widetilde{X}$ and $\widetilde{Y}$ of $X$ and $Y$, respectively, and an extension $\widetilde{f}\colon\widetilde{X}\to\widetilde{Y}$ of $f$ such that $\widetilde{f}$ is closed, $\dim \widetilde{f}=\dim f$, $\dim\widetilde{X}=\dim X$ and $\dim\widetilde{Y}=\dim Y$. In the present note we extend this result for extensional dimension with respect to countable $CW$-complexes. We also establish an analogue (see Theorem 3.1) of a result of Katetov \cite{mk} characterizing the dimension $\dim$ of metrizable spaces in terms of uniformly 0-dimensional maps into finite-dimensional cubes.
 
Recall that $e-dim X\leq K$ if and only if  every continuous map $g\colon A\to K$, where $A\subset X$ is closed, can be extended to a map $\bar{g}\colon X\to K$, see  \cite{d:95}. For a map $f\colon X\to Y$ we write
$e-dim f\leq K$ provided $e-dim f^{-1}(y)\leq K$ for every $y\in Y$.  Unless indicated otherwise, all spaces are assumed to be metrizable and all maps continuous. By a $CW$-complex we always mean a countable $CW$-complex.


\section{Completion of maps}

We begin with the following lemma.

\begin{lem} Let  $f\colon X\to Y$ be a perfect map between metrizable spaces and $K$ a $CW$-complex. Then  
$B_K=\{y\in Y: \e f^{-1}(y)\leq K\}$ is a $G_{\delta}$-subset of $Y$.
\end{lem}

\begin{proof} 
By \cite{bp:98}, there exists a map $g$ from $X$ into the Hilbert cube $Q$ such that $f\times g\colon X\to Y\times Q$ is an embedding. Let $\{W_i\}_{i\in\N}$ be  a countable finitely-additive base for $Q$. 
For every $i$ we choose a sequence of mappings $h_{ij}\colon\overline{W_i}\to K$, representing all the homotopy classes of mappings from $\overline{W_i}$ to $K$ (this is possible because $K$ is a countable $CW$-complex and all $\overline{W_i}$ are metrizable compacta ).  For any $i$, $j$ let $U_{ij}$ be the set of all  $y\in Y$ having the following property:

\smallskip\noindent
the map $h_{ij}\circ g\colon g^{-1}(\overline{W_{i}})\to K$ can be continuously extended to a map over the set
$g^{-1}(\overline{W_{i}})\cup  f^{-1}(y)$.

\smallskip\noindent 
Let show that every $U_{ij}$ is open in $Y$. Indeed, if $y_0\in U_{ij}$, then there exists a map $h\colon g^{-1}(\overline{W_{i}})\cup  f^{-1}(y_0)\to K$ extending $h_{ij}\circ g$. Since $K$ is an absolute extensor for metrizable spaces, we can extend $h$ to a map $\overline{h}\colon V\to K$, where $V\subset X$ is open and contains  $g^{-1}(\overline{W_{i}})\cup  f^{-1}(y_0)$. Because $f$ is 
closed, there exists a neighborhood $G$ of $y_0$ in $Y$ with $f^{-1}(G)\subset V$. Then, for every $y\in G$, the restriction of $\overline{h}$ on $g^{-1}(\overline{W_{i}})\cup  f^{-1}(y)$ is an extension of $h_{ij}\circ g$. Hence, $G\subset U_{ij}$.

\smallskip\noindent
It is clear that $B_K$ is contained in every $U_{ij}$. It remains only to show that $\cap_{i,j=1}^{\infty}U_{ij}\subset B_K$. Take
$y\in \cap_{i,j=1}^{\infty}U_{ij}$ and a map $h\colon A\to K$, where $A$ is a closed subset of $f^{-1}(y)$. 
Because the map $g_y=g|f^{-1}(y)$ is a homeomorphism, $h^{'}=h\circ g_y^{-1}\colon g(A)\to K$ is well defined. Next, extend $h^{'}$ to a map from a neighborhood  $W$ of $g(A)$ in $Q$ (recall that $f^{-1}(y)$ is compact, so $g(A)\subset Q$ is closed) into $K$ and find $W_k$ with $g(A)\subset W_k\subset\overline{W_k}\subset W$. Therefore, there exists a map
$h^{''}\colon\overline{W_k}\to K$ extending $h^{'}$.  Then $h^{''}$ is homotopy equivalent to some $h_{kj}$, so are 
$h^{''}\circ g$ and $h_{kj}\circ g$ (considered as maps from $g^{-1}(\overline{W_{k}})$ into $K$).  Since $y\in U_{kj}$, $h_{kj}\circ g$ can be extended to a map from $g^{-1}(\overline{W_{k}})\cup f^{-1}(y)$ into $K$. Then, by the Homotopy Extension Theorem, there exists a map $\bar{h}\colon g^{-1}(\overline{W_{k}})\cup f^{-1}(y)\to K$ extending $h^{''}\circ g$. Obviously, $\bar{h}|f^{-1}(y)$ extends  $h$. Hence, 
$\e f^{-1}(y)\leq K$. 
\end{proof}

The next lemma, though not explicitely stated in this form, was actually proved by Levin \cite{ml1}.

\begin{lem}
Let $X$ be a subset of the metrizable space $Y$ with $\e X\leq K$ for some $CW$-complex $K$. Then there exists a $G_{\delta}$-subset $\widetilde{X}$ of $Y$ containing $X$ such that $\e\widetilde{X}\leq K$.
\end{lem} 

\begin{thm} Let $f\colon X\to Y$ be a closed surjective map between metrizable spaces such that $\e f\leq K$, $\e X\leq L_X$ and $\e Y\leq L_Y$, where $K$, $L_X$ and $L_Y$ are $CW$-complexes. Then there exist completions $\widetilde{X}$ and $\widetilde{Y}$ of $X$ and $Y$, respectively, and a closed surjection $\widetilde{f}\colon\widetilde{X}\to\widetilde{Y}$ extending $f$ with $\e\widetilde{f}\leq K$, $\e\widetilde{X}\leq L_X$ and $\e\widetilde{Y}\leq L_Y$. 
\end{thm}

\begin{proof}  
Since $f$ is closed, $Fr f^{-1}(y)=\emptyset$ if and only if $y$ is a discrete point in $Y$, where $Fr f^{-1}(y)$ denotes the boundary of $f^{-1}(y)$ in $X$. On the other hand, it is easily seen that the validity of the theorem for any metrizable $Y$ without discrete points implies its validity for any metrizable $Y$. Therefore, we can assume that $Y$ doesn't have any discrete points, or equivalently, $Fr f^{-1}(y)\neq\emptyset$ for every $y\in Y$. According to the classical result of Va\v{i}nste\v{i}n \cite{v} (see also \cite{jk}), there are completions $X_1$ and $Y_1$ of $X$ and $Y$, respectively, and a closed surjection $f_1\colon X_1\to Y_1$ which extends $f$. For any $y\in Y_1$ we denote by $Fr f_1^{-1}(y)$ the boundary of $f_1^{-1}(y)$ in $X_1$. Then, the following two facts occur:

\medskip\noindent
(1)\hbox{}~~~~~~$Fr f^{-1}(y)$ coincides with $Fr f_1^{-1}(y)$ provided $y\in Y$;

\smallskip\noindent
(2)\hbox{}~~~~~~$f_1^{-1}(y)=Fr f_1^{-1}(y)$ provided $y\in Y_1\backslash Y$. 

\medskip\noindent
Therefore, $Fr f_1^{-1}(y)\neq\emptyset$ for all $y\in Y_1$.  Moreover, $f_1|H\colon H\to Y_1$ is a perfect surjection (see \cite{v}), where $H=\bigcup\{Fr f_1^{-1}(y): y\in Y_1\}$. Obviously, $H$ is closed in $X_1$, so $\e (H\cap X)\leq L_X$. Then, by Lemma 2.2, there exists a $G_{\delta}$-subset $P$ of $H$ with $H\cap X\subset P$ and 

\medskip\noindent
(3)\hbox{}~~~~~~$\e P\leq L_X$. 

\medskip\noindent
It follows from (1) that $(f_1|H)^{-1}(Y)\subset P$. Therefore, $f_1(H\backslash P)$ does not meet $Y$. Since  $f_1|H$ is a closed surjection onto $Y_1$ and $H\backslash P$ is $F_{\sigma}$ in $H$, $f_1(H\backslash P)$ is $F_{\sigma}$ in $Y_1$.
So, $Y_2=Y_1\backslash f_1(H\backslash P)$ is a $G_{\delta}$-set  in $Y_1$ containing $Y$ and such that 

\medskip\noindent
(4)\hbox{}~~~~~~$(f_1|H)^{-1}(Y_2)\subset P$.

\medskip\noindent   
Condition (1) also implies that every fiber $(f_1|H)^{-1}(y)$ is of extensional dimension $\leq K$ provided $y\in Y$. Hence, applying Lemma 2.1 and then Lemma 2.2, we can find a $G_{\delta}$-subset $\widetilde{Y}$ of $Y_2$ such that $\e \widetilde{Y}\leq L_Y$ and

\medskip\noindent
(5)\hbox{}~~~~~~$\e Fr f_1^{-1}(y)\leq K$ for all $y\in\widetilde{Y}$.

\medskip\noindent
Consider the set $W=X_1\backslash H$. It is open in $X_1$, so $W\cap X$ is open in $X$. Moreover, $f^{-1}(y)\cap W$ is the interior of $f^{-1}(y)$ in $X$,  $y\in Y$. Therefore, $\e \big(f^{-1}(y)\cap W\big)\leq K$ for every $y\in Y$. Consequently, $\e \big(W\cap X\big)\leq K$. On the other hand, $W\cap X$ is a subset of $X$, so $\e \big(W\cap X\big)\leq L_X$. Since the property of metrizable spaces to have extensional dimension less than or equal to a given countable $CW$-complex is hereditary (see, for example \cite{ch}), we can apply Lemma 2.2 twice to obtain a $G_{\delta}$-subset $U$ of $W$ which contains  $W\cap X$ such that 

\medskip\noindent
(6)\hbox{}~~~~~~$\e U\leq K$ and $\e U\leq L_X$. 

\medskip\noindent
Finally, let $\widetilde{X}=f_1^{-1}(\widetilde{Y})\cap (U\cup P)$ and $\widetilde{f}=f_1|\widetilde{X}$. Obviously, $\widetilde{X}\cap U$ and $\widetilde{X}\cap P$ are disjoint, respectively, open and closed subsets of 
$\widetilde{X}$. Since $\e \big(\widetilde{X}\cap U\big)\leq\e U\leq L_X$ and 
$\e \big(\widetilde{X}\cap P\big)\leq\e P\leq L_X$, $\widetilde{X}$ can be represented as the union of countable many its closed subsets $F_i$ with $\e F_i\leq L_X$ for each $i$. Then, by the countable sum theorem, $\e\widetilde{X}\leq L_X$. It follows from our construction that $\widetilde{f}$ maps $\widetilde{X}$ onto $\widetilde{Y}$ and each $\widetilde{f}^{-1}(y)$, $y\in\widetilde{Y}$, is the union of the disjoint sets $Fr f_1^{-1}(y)$ and $\widetilde{f}^{-1}(y)\cap U$ which are, respectively, closed and open in $\widetilde{f}^{-1}(y)$. 
By (5) and (6), both $Fr f_1^{-1}(y)$ and $\widetilde{f}^{-1}(y)\cap U$ are of extensional dimension $\leq K$. Hence, $\e f_1^{-1}(y)\leq K$ for each $y\in\widetilde{Y}$.

\medskip\noindent
It only remains to show that $\widetilde{f}$ is a closed map. To this end, let $A\subset\widetilde{X}$ be closed and $y_n=\widetilde{f}(x_n)$ converges to $y_0$, where $\{x_n\}$ is a sequence of points from $A$. Suppose that $y_0\not\in\widetilde{f}(A)$. Then, by (1), (2) and (4), $Fr f_1^{-1}(y_0)\subset\widetilde{X}$ and it does not meet $A$ (as a subset of $\widetilde{f}^{-1}(y_0)$). Being compact $Fr f_1^{-1}(y_0)$ is closed in $\widetilde{X}$. Consequently, there is an open $V\subset X_1$ containing $Fr f_1^{-1}(y_0)$ such that $V\cap A=\emptyset$. 
Let $V_1$ be the union of $V$ and the interior of $f_1^{-1}(y_0)$ in $X_1$. Obviously, $V_1$ is open in $X_1$, contains $f_1^{-1}(y_0)$ and does not meet $A$. 
Since $f_1$ is a closed map, there exists a neighborhood $O(y_0)$ of $y_0$ in $Y_1$ such that $f_1^{-1}(y)\subset V_1$ for all $y\in O(y_0)$. Therefore, $f_1^{-1}(y_m)\subset V_1$ for some $m$. The last inclusion implies $x_m\in V_1\cap A$, which is a contradiction. Therefore, $y_0\in\widetilde{f}(A)$, i.e. $\widetilde{f}$ is closed.
\end{proof}


\section{$\sigma$-uniformly 0-dimensional maps}

A map $f\colon X\to Y$ is called uniformly 0-dimensional \cite{mk} if there exists a metric on $X$ generating its topology such that for every $\epsilon>0$ every point of $f(X)$ has a neighborhood $U$ in $Y$ with $f^{-1}(U)$ being the union of disjoint open subsets of $X$ each of diameter $<\epsilon$.  Uniformly 0-dimensional maps are called in \cite{ap:73} completely 0-dimensional.  
It is well known, that if $f\colon X\to Y$ is uniformly 0-dimensional and $\dim Y\leq n$, then $\dim X\leq n$  (see, for example, \cite{mk}, \cite{ap:73} or \cite{ml1}).

We say that a map $g\colon X\to Y$ is $\sigma$-uniformly 0-dimensional if $X$  can be represented as the union of countably many of its closed subsets $X_i$ such that each restriction $g|X_i$ is uniformly 0-dimensional.  Katetov \cite{mk} (see also \cite{jn}) proved that a space $X$  is at most $n$-dimensional if and only if for each metrization of $X$ there exists a uniformly 0-dimensional map of $X$ into $\uin^n$. Moreover, the space $C(X,\uin^n)$ with the uniform convergence topology contains a dense $G_{\delta}$-subset  consisting of uniformly 0-dimensional maps. The next theorem can be considered as a parametric version of Katetov's result (see \cite{re:95} for the definition of $C$-spaces).

\begin{thm} Let $f\colon X\to Y$ be a closed map of metrizable spaces with $Y$ being a $C$-space. Then $\dim f\leq n$ if and only if 
there exists a map $g\colon X\to \uin^n$ such that $f\times g$ is $\sigma$-uniformly 0-dimensional. Moreover, if $\dim f\leq n$, then the set of all such maps $g\in C(X,\uin^n)$ is dense in $C(X,\uin^n)$ with respect to the uniform convergence topology generated by the Euclidean metric on $\uin^n$.
\end{thm}

\begin{proof} All function spaces in this proof are equipped with the uniform convergence topology.  

Suppose that $\dim f\leq n$.
We represent $X$ as the union $X=X_0\cup (X\backslash X_0)$ such that
$X_0$ is closed in $X$, $f_0=f|X_0$ is a perfect map and $\dim (X\backslash X_0)\leq n$.  Let $X\backslash X_0=\bigcup_{k=1}^{\infty}X_k$ such that each $X_k$ is closed in $X$.  Since $f_0\colon X_0\to Y$ is perfect, the set $C_0$ of all $g\colon X\to\uin^n$ with $(f\times g)|X_0$ being 0-dimensional is dense in  $C(X,\uin^n)$ (see for example, \cite[Theorem  1.3]{tv:02}).  It is easily seen that every perfect 0-dimensional map between metric spaces is uniformly 0-dimensional. Hence, all restrictions $(f\times g)|X_0$, $g\in C_0$, are uniformly 0-dimensional.  For every $g\in C_0$ let
$H(g)=\{h\in C(X,\uin^n): h|X_0=g|X_0\}$.  Each $H(g)$ is  closed in $C(X,\uin^n)$ and $C_0=\cup\{H(g): g\in C_0\}$.
We also define the maps $p_k\colon C(X,\uin^n)\to C(X_k,\uin^n)$ by $p_k(h)=h|X_k$,  $k=1,2,..$, and let $p_{k,g}\colon H(g)\to C(X_k,\uin^n)$ denote the restriction $p_k|H(g)$ for any $k\in\N$ and $g\in C_0$.  Using that $X_0$ and each $X_k$ are disjoint closed sets in $X$, we can show that every $p_{k,g}$ is open and surjective.  According to the Katetov result  \cite{mk}, there exists a dense and $G_{\delta}$-subset $C_k$ of $C(X_k,\uin^n)$ consisting of uniformly 0-dimensional maps, $k=1,2,..$. Consequently, for any $g\in C_0$, the sets $H_k(g)=p_{k,g}^{-1}(C_k)$ are dense and $G_{\delta}$ in $H(g)$. Since $H(g)$ has the Baire property (as a closed subset of $C(X,\uin^n)$),  $M(g)=\bigcap_{k=1}^{\infty}H_k(g)$ is also dense and $G_{\delta}$ in $H(g)$.
Then $M=\cup\{M(g): g\in C_0\}$ is dense in $C(X,\uin^n)$. Moreover, it follows from the construction that, for any $g\in M$, the restrictions $(f\times g)|X_k$ are uniformly 0-dimensional, $k=0,1,2,..$. Therefore, $M$ consists of $\sigma$-uniformly 0-dimensional maps. 

To prove the other implication of Theorem 3.1, assume that there exists $g\colon X\to\uin^n$ such that the map $f\times g\colon X\to Y\times\uin^n$ is $\sigma$-uniformly 0-dimensional. Therefore, $X$ can be represented as the union of countably many of its closed subsets $A_i$ such that each $(f\times g)|A_i$ is uniformly 0-dimensional.  The last condition implies that, for any $y\in Y$ and $i$ the map $g|(f^{-1}(y)\cap A_i)\colon f^{-1}(y)\cap A_i\to \uin^n$ is uniformly 0-dimensional. Hence, $\dim\big(f^{-1}(y)\cap A_i\big)\leq n$. Since $f^{-1}(y)=\cup_{i=1}^{\infty}f^{-1}(y)\cap A_i$,  by the countable sum theorem, $\dim f^{-1}(y)\leq n$ for each $y\in Y$. So, $\dim f\leq n$.
\end{proof}


\bigskip

\end{document}